\input psfig.tex
\language=3

\def\nobf{\noindent\bf}

\tolerance=5000
\magnification=\magstep1

\parskip=4pt plus 6pt

\font\bigbf=cmb10 scaled \magstep 2
 3

\input amssym.def
\input amssym.tex

\def\CC{{\Bbb C}}

\def\RR{{\Bbb R}}
\def\SS{{\Bbb S}}

\def\ZZ{{\Bbb Z}}

\def\fullsquare{{\vrule height6pt width6pt depth0pt}}

\def\qed{\unskip\nobreak
\hfil\penalty50\hskip1.75em\null\nobreak\hfil\fullsquare
{\parfillskip=0pt \finalhyphendemerits=0 \par}}

%

\newcount\m \newcount\n \newcount\p \newdimen\dim

\def\today{\ifcase\month\or
  January\or February\or March\or April\or May\or June\or
  July\or August\or September\or October\or November\or December\fi
  \space\number\day, \number\year}

\def\hoje{\number\day\ de
  \ifcase\month\or
  janeiro\or fevereiro\or mar{\c c}o\or abril\or maio\or junho\or
  julho\or agosto\or setembro\or outubro\or novembro\or dezembro\fi
  \ de \number\year}

\def\japtoday{\number\year.\twodigits\month.\twodigits\day}

\def\hours{\n=\time \divide\n 60
  \m=-\n \multiply\m 60 \advance\m \time
  \twodigits\n:\twodigits\m}

\def\twodigits#1{\ifnum #1<10 0\fi \number#1}

\def\datepages{\footline={{\hss\tenrm\folio\hss\hbox to 0pt{
\hidewidth\fiverm\japtoday}}}}

\font\eightrm cmr8

\def\aa{{\bf a}}
\def\kk{{\bf k}}
\def\pp{{\bf p}}
\def\cG{{\cal G}}
\def\cH{{\cal H}}
\def\bcG{\bar{\cal G}}
\def\bcH{\bar{\cal H}}
\def\Zz{{\ZZ/(2)}}
\def\AA{{\bf A}}
\def\DD{{\bf D}}

\def\Cq{{\CC_q}}
\def\Sq{{\SS_q}}
\def\qq{{\bf q}}

\centerline{\bigbf Singular polynomials of}
\vskip5pt
\centerline{\bigbf generalized Kasteleyn matrices}

\medskip

\centerline{\it Nicolau C. Saldanha%
\footnote{}{\eightrm Research supported by CNPq (Brazil) and ENS-Lyon}}

\bigskip

{

\narrower

\noindent{\bf Abstract: }
Kasteleyn counted the number of domino tilings of a rectangle
by considering a mutation of the adjacency matrix:
a {\sl Kasteleyn matrix} $K$.
In this paper we present a generalization of Kasteleyn matrices
and a combinatorial interpretation
for the coefficients of the characteristic polynomial of $KK^\ast$
(which we call the {\sl singular polynomial}),
where $K$ is a generalized Kasteleyn matrix for a planar bipartite graph.
We also present a $q$-version of these ideas and a few results
concerning tilings of special regions such as rectangles.

}

\bigskip

{\nobf Introduction}

\smallskip

Kasteleyn ([K]) counted the number of domino tilings of a rectangle
by considering a mutation of the adjacency matrix,
since then known as a {\sl Kasteleyn matrix} ([LL], [ST]).
Given a planar bipartite graph $\cG$
there are several Kasteleyn matrices $K$ for $\cG$
but, as has been shown independently by David Wilson and Horst Sachs,
the singular values of $K$ or, equivalently,
the eigenvalues of $KK^\ast$, are independent of the choice of $K$.
Following a question posed by James Propp ([P]),
we search for a combinatorial interpretation for these numbers.

In section 1 we introduce generalized Kasteleyn matrices
for planar bipartite graphs
and present a combinatorial interpretation for the determinant
of such matrices $A$ in terms of counting matchings.
In section 2 we address the main issue of understanding
what the coefficients of the characteristic polynomial of $AA^\ast$ represent,
and then, in section 3, we consider the special case of Kasteleyn matrices.
In section 4 we present the $q$-analogs of these ideas.
In section 5 we take a look at rectangles in the plane.
Finally, in section 6 we present a few other small examples.
We find the language of homology theory helpful and use it throughout the paper.
We thank James Propp, Richard Kenyon and Horst Sachs for helpful
conversations and emails.

\bigbreak

{\nobf 1. Generalized Kasteleyn matrices and their determinants}

\nobreak\smallskip\nobreak

Let $\cG$ be a planar bipartite graph with $n$ white vertices
and $n'$ black vertices.
We number the white (black) vertices $1, 2, \ldots, n$ ($1, 2, \ldots, n'$).
A {\sl generalized Kasteleyn matrix} for $\cG$
is an $n \times n'$ complex matrix $A$ such that
$$|a_{ij}| = \cases{1,&if the $i$-th white vertex
and the $j$-th black vertex are adjacent,\cr
0& otherwise.\cr}$$
Such matrices are conveniently represented
by labeling the edges of $\cG$ with complex numbers of norm 1.

We may identify a generalized Kasteleyn matrix $A$
with a cocomplex $\AA \in C^1(\cG, \SS^1)$
by making the convention that $a_{ij}$ indicates the value of $\AA(e_{ij})$,
$e_{ij}$ being the oriented edge
going from the $j$-th black to the $i$-th white vertex.
A notational confusion must be avoided here:
the complex numbers of norm 1 form a multiplicative group
but the coefficients for homology or cohomology should be additive groups.
Thus, from now on, the symbol $\SS^1$ shall denote the additive group $\RR/\ZZ$
and we denote the exponential $x \mapsto \exp(2 \pi i x)$
by $\eta: \SS^1 \to \CC$.
In particular, we write $a_{ij} = \eta(\AA(e_{ij}))$.
Since $C^2(\cG, \SS^1) = 0$,
any cocomplex $\AA$ is automatically closed and
a generalized Kasteleyn matrix $A$ defines an element $\aa \in H^1(\cG,\SS^1)$.

There is a natural inclusion $\Zz \subseteq \SS^1$;
this defines $\eta: \Zz \to \CC$ with $\eta(m) = (-1)^m$.
We also obtain induced inclusions
$C^1(\cG, \Zz) \subseteq C^1(\cG, \SS^1)$ and
$H^1(\cG, \Zz) \subseteq H^1(\cG, \SS^1)$.
For a generalized Kasteleyn matrix $A$,
$\AA \in C^1(\cG, \Zz)$ if and only if $A$ is a real matrix.

{\nobf Lemma 1.1: }
{\sl There is a unique element $\kk \in H^1(\cG, \Zz)$
such that for any cycle $C$,
$$\kk(C) \equiv m + l + 1 \pmod 2 \eqno{(1)}$$
where $m$ is the number of vertices in the interior of $C$
and $2l$ is the length of $C$. }

In the statement above, the word `cycle' is used in the sense of graph theory:
$C$ is a simple closed curve in the plane
composed of edges and vertices of $\cG$
and the interior of $C$ is well defined by Jordan's theorem.
Of course, graph theory cycles define homology cycles
(i.e., closed elements of $C_1(\cG,\ZZ)$) but the converse is not always true;
any homology cycle may nevertheless be written a linear combination
of graph theory cycles.

{\nobf Proof: }
Uniqueness is obvious since the above equation gives the value of $\kk$
computed against any cycle and thus, by linearity, against any element
of $H_1(\cG, \ZZ)$
(recall that $H^1(\cG, \Zz) = {\rm Hom}(H_1(\cG, \ZZ), \Zz)$).

In order to construct $\kk$, we first notice that
$H^1(\cG,\Zz) = (\Zz)^h$, where $h$ is the number of {\sl holes}
(bounded connected components of the complement) of $\cG$:
in this identification, the coordinates of $\aa$
corresponding to a given hole is $\aa(C)$,
$C$ being the outer boundary of the said hole.
It is then clear that there exists a unique element of $H^1(\cG,\Zz)$,
which we call $\kk$, satisfying equation $(1)$
for all such $C$.

Figure 1.1 illustrates a minor complication which has to be kept in mind:
the boundary of a hole is not always a cycle in the sense of graph theory.
There should be no confusion, however:
in (a) we have $l = 2$, $m = 1$ for the only hole
and in (b) we have $l = 2$, $m = 0$ and $l = 2$, $m = 4$
for the smaller and bigger hole, respectively.

\midinsert \smallskip

\centerline{%
\vbox{\hsize=25mm%
\centerline{\psfig{file=diag0a.eps}}
\bigskip \centerline{\eightrm (a)}}\qquad%
\vbox{\hsize=25mm%
\centerline{\psfig{file=diag0b.eps}}
\bigskip \centerline{\eightrm (b)}}\qquad%
\vbox{\hsize=25mm%
\centerline{\psfig{file=diag0c.eps}}
\bigskip \centerline{\eightrm (c)}}}

\smallskip

\centerline{\eightrm Figure 1.1}

\smallskip \endinsert

Let $C$ be an arbitrary cycle: we prove that equation (1) holds for $C$.
The interior of $C$ minus $\cG$ is a union of holes.
If we discard the holes which are completely surrounded by other holes in $C$
and consider the outer boundaries $C_1, \ldots, C_k$ of the remaining ones,
we have $\kk(C) = \sum_i \kk(C_i)$.
Since equation (1) holds for each $C_i$ we have
$$\kk(C) \equiv k + \sum_i l_i + \sum_i m_i \pmod 2,$$
where $l_i$ and $m_i$ correspond to $C_i$.
Notice that $L = l + \sum_i l_i$,
$L$ denoting the number of edges of some $C_i$
on $C$ or in the interior of $C$,
and $M = 2l + m - \sum_i m_i$,
$M$ denoting the number of vertices of some $C_i$
on $C$ or in the interior of $C$.
Finally, by Euler characteristic, $k - L + M = 1$
and we have equation $(1)$ for $C$, proving our lemma.
\qed

We call $\kk \in H^1(\cG, \Zz)$ as defined in the previous lemma
the {\sl Kasteleyn class};
when the graph $\cG$ is not clear from the context, we write $\kk_\cG$.
The definition of $\kk$ involves $m$ and thus appears to depend
on the way $\cG$ is drawn in the plane.
Indeed, examples (a) and (c) in Figure 1.1 represent equivalent graphs,
but the Kasteleyn classes are different;
solid lines stand for a label 1 and dashed lines stand for a label $-1$.
A {\sl Kasteleyn matrix} is a generalized Kasteleyn matrix
corresponding to the Kasteleyn class.

\bigskip

We restrict ourselves for the rest of this section to the case $n = n'$
in order to explore the relationship between matchings
and the determinant of square generalized Kasteleyn matrices $A$.

It is natural to interpret a matching of $\cG$
as the sum of its edges oriented from black to white
and thus as an element of $C_1(\cG, \ZZ)$.
The boundary of any matching always equals
the sum of all white vertices minus the sum of all black vertices;
thus, the difference of two matchings of $\cG$ is closed
and may be identified with an element of $H_1(\cG, \ZZ)$.
Notice furthermore that the difference of two matchings may be written
in a unique way as a sum of disjoint graph theory cycles.

For $\aa \in H^1(\cG, \SS^1)$ and two matchings $\mu_1$ and $\mu_2$,
$\eta(\aa(\mu_2 - \mu_1))$ is a complex number of absolute value 1.
In particular, if $\aa \in H^1(\cG, \Zz)$ then
$\eta(\aa(\mu_2 - \mu_1))$ is 1 or $-1$.
We then say that $\mu_1$ and $\mu_2$
have the same {\sl $\aa$-parity} if $\eta(\aa(\mu_2 - \mu_1)) = 1$;
$\aa$-parity splits the set of matchings into two equivalence classes
(occasionally one of these classes may turn out to be empty).

{\nobf Lemma 1.2: }
{\sl For any planar graph $\cG$, for any $\aa \in H^1(\cG, \SS^1)$ 
and for any given matching $\mu_0$ we have
$$\sum_{\mu_1,\mu_2} \eta(\aa(\mu_2 - \mu_1)) =
\left|{\sum_\mu \eta(\aa(\mu - \mu_0))}\right|^2,$$
where $\mu_1, \mu_2$ and $\mu$ range over all matchings.}

{\nobf Proof: }
We may write the right hand side as
$$
\left({\sum_{\mu_2} \eta(\aa(\mu_2 - \mu_0))}\right)
\left({\sum_{\mu_1} \overline{\eta(\aa(\mu_1 - \mu_0))}}\right) =
\left({\sum_{\mu_2} \eta(\aa(\mu_2 - \mu_0))}\right)
\left({\sum_{\mu_1} \eta(\aa(\mu_0 - \mu_1))}\right)
$$
and distribute to get the left hand side.
\qed

Define
$$\delta(\aa,\cG) = \sum_{\mu_1,\mu_2} \eta(\aa(\mu_1 - \mu_2)),$$
where $\mu_1$ and $\mu_2$ range over all matchings;
if $\cG$ admits no matchings we define $\delta(\aa,\cG) = 0$.
As an example, $\delta(0,\cG)$
is the square of the number of matchings of $\cG$.
Also, for $\aa \in H^1(\cG, \Zz)$,
$\delta(\aa, \cG)$ is the square of the difference
between the number of matchings in each $\aa$-parity equivalence class.

A matching may also be though of as a bijection
from the set of white vertices to the set of black vertices.
Thus, if $\mu_1$ and $\mu_2$ are matchings then $\mu_1^{-1} \circ \mu_2$
is a permutation of the set of white vertices:
we say that these two matchings have the same {\sl permutation parity}
if and only if the this permutation is even.

{\nobf Lemma 1.3: }
{\sl Two matchings have the same permutation parity if and only if
they have the same $\kk$-parity, $\kk$ being the Kasteleyn class. }

{\nobf Proof: }
Let $\mu_1$ and $\mu_2$ be two matchings
and write $\mu_1 - \mu_2$ as a sum of disjoint cycles
$C_1,\ldots,C_N$ of lengths $2l_1,\ldots,2l_N$.
The interior and exterior of any of these cycles is matchable,
thus $m_1,\ldots,m_N$ as in Lemma 1.1 are all even.
From equation $(1)$, $\kk(C_i) \equiv l_i + 1 \pmod 2$ and thus
$\kk(\mu_0 - \mu_1) \equiv \sum (l_i + 1) \pmod 2$.

The permutation $\mu_0^{-1} \circ \mu_1$ can be written as a product
of $N$ cycles (in the permutation sense) corresponding to $C_1,\ldots,C_N$
with lengths $l_1,\ldots,l_N$ and the parity of the permutation
$\mu_1^{-1} \circ \mu_2$ is thus $\sum (l_i + 1)$.
This proves our claim.
\qed

Notice that permutation parity, unlike the Kasteleyn class,
does not depend on how $\cG$ is drawn in the plane.
A corollary of the previous lemma is thus that if differences of matchings
generate $H_1(\cG,\ZZ)$ then the Kasteleyn class of $\cG$
is the same for all planar embeddings.

{\nobf Lemma 1.4: }
{\sl For any generalized Kasteleyn matrix $A$ we have
$|\det(A)|^2 = \delta(\aa + \kk,\cG)$. }

{\nobf Proof: }
Each non-zero monomial in the expansion of $\det(A)$
corresponds to a matching.
Thus, each matching $\mu$ contributes with
a complex number of absolute value 1 to $\det(A)$.
The expression $\eta(\aa(\mu - \mu_0))$ obtains,
up to a fixed multiplicative constant of absolute value 1,
the product of the corresponding elements of $A$.
From Lemma 1.2, $\kk$-parity is permutation parity, i.e.,
gives the sign of the monomial in the definition of the determinant.
Thus, the contribution of $\mu$ to $\det(A)$ is,
again up to a fixed multiplicative constant of absolute value 1,
$\eta((\aa + \kk)(\mu - \mu_0))$, proving our lemma.
\qed

As a special case, if $K$ is a Kasteleyn matrix,
$|\det(K)|$ is the number of matchings of $\cG$:
this is Kasteleyn's original motivation.

\bigbreak

{\nobf 2. Singular polynomials of generalized Kasteleyn matrices}

\nobreak\smallskip\nobreak

Having provided an interpretation for $|\det(A)|$ when $A$ is square,
it is natural to ask about other functions of $A$,
specially if $A$ is not square.
We should not expect natural interpretations for the argument of $\det(A)$
since it depends on the way we assign labels to vertices.
Also, a few simple experiments will show that
the spectrum of $A$ (even if $A$ is square) is not a function of $\aa$.
The following lemma tells us what functions of $A$
are determined by $\aa$.

{\nobf Lemma 2.1: }
{\sl Let $A$ be a generalized Kasteleyn matrix for $\cG$
and let $\aa$ be the corresponding element of $H^1(\cG, \SS^1)$.
Then the generalized Kasteleyn matrices for $\cG$
also corresponding to $\aa$ are precisely the matrices of the form
$D_1 A D_2$ where $D_1$ and $D_2$ are unitary diagonal matrices.
Furthermore, if $\cG$ is connected, $D_1 A D_2 = D'_1 A D'_2$
if and only if there exists a complex number $z$ of absolute value 1 with
$D_1 = z D'_1$, $D_2 = z^{-1} D'_2$. }

It is possible to give a more elementary proof,
but following the spirit of the rest of this paper
we phrase the proof in homological language.

{\nobf Proof: }
As we saw in Section 1, generalized Kasteleyn matrices
correspond to 1-cocomplexes in $C^1(\cG,\SS^1)$;
two such 1-cocomplexes $\AA$ and $\AA'$
induce the same element of $H^1(\cG, \SS^1)$
if and only if their difference is the coboundary of a 0-cocomplex.
A 0-cocomplex $\DD$ is a function assigning
an element of $\SS^1$ to each vertex;
the $\eta$'s of these elements 
may conveniently be arranged in a pair of unitary diagonal matrices,
$D_w$ for the white and $D_b$ for the black vertices.
It is a simple translating process to verify that
the cocomplex $\AA + d(\DD)$ corresponds
to the generalized Kasteleyn matrix $D_w A D_b^{-1}$,
thus proving our first claim.
The uniqueness of $D_1$ and $D_2$ up to a constant multiplicative factor
corresponds to the fact that the only closed 0-cocomplexes are the constants,
i.e., that $H^0(\cG, \SS^1) = \SS^1$ if $\cG$ is connected.
\qed

We recall that for any complex $n \times n'$ matrix $B$,
there are unitary matrices $U_1$ and $U_2$ such that
$S = U_1 B U_2$ is a real diagonal matrix
with non-increasing non-negative diagonal entries;
$S$ is well-defined given $B$ and its diagonal entries
(i.e., the $s_{ii}$ entries of $S$, even if $S$ is not square)
are called the {\sl singular values} of $B$.
The rows of $U_1$ (resp., columns of $U_2$)
are called the {\sl left} (resp., {\sl right}) {\sl singular vectors} of $B$.
It is easy to see that the singular values and
left (resp. right) singular vectors of $B$
are the non-negative square roots of the eigenvalues 
and the eigenvectors of $BB^\ast$ (resp., $B^\ast B$).
Inspired in these classical notions, we call the characteristic polynomial
of $BB^\ast$ the {\sl singular polynomial} of $B$:
its roots are the squares of the singular values of $B$.
Also, the singular polynomials of $B$ and $U_1 B U_2$ are equal
and the singular polynomials of $B$ and $B^\ast$
differ by a factor of $t^{n-n'}$.

It follows from Lemma 2.1 and the remarks in the previous paragraph
that the singular polynomial of $A$ is determined by $\aa$: we call it $P_\aa$.
The singular values of $A$ and, if the singular values are simple,
the absolute values of the coordinates of the singular vectors
(up to a constant factor) are also determined by $\aa$.
We shall now present what we find to be a reasonably natural interpretation
for the coefficients of $P_\aa$.
While these numbers determine the singular values
the question remains whether a nice interpretation exists
for the actual singular values and vectors.

\smallskip

Let $\cH \subseteq \cG$ be a balanced subgraph of $\cG$:
the inclusion induces a map
$\pi_{\cG,\cH}: H^1(\cG, \SS^1) \to H^1(\cH, \SS^1)$.
More concretely, if $\aa$ corresponds to a generalized Kasteleyn matrix $A$
then $\pi_{\cG,\cH}(\aa)$ corresponds to the submatrix of $A$
obtained by picking only the elements for which both row and column
correspond to elements of $\cH$.
The simplest interpretation is probably in terms of labels for edges:
just keep the old labels.
When this causes no confusion, we write $\aa$ instead of $\pi_{\cG,\cH}(\aa)$:
for instance, we write $\delta(\aa, \cH)$ instead of the more correct
but cumbersome $\delta(\pi_{\cG,\cH}(\aa), \cH)$.

{\nobf Theorem 2.2: }
{\sl Let $A$ be a generalized Kasteleyn matrix and let
$P_\aa(t) = t^n + a_1 t^{n-1} + \cdots + a_{n-1} t + a_n$
be the singular polynomial of $A$.
Then 
$$a_m = (-1)^m \sum_{|\cH| = 2m} \delta(\aa + \kk_\cH, \cH)\eqno{(2)}$$
where $\cH$ ranges over all balanced subgraphs with $2m$ vertices.}

Notice that for $m = n$ equation (2) is equivalent to Lemma 1.4.
For $m = 1$, we get the simple remark that
$|a_1|$ is the number of edges of $\cG$, regardless of $\aa$.
An interpretation for $a_2$ is already subtler:
each subgraph with two white and two black vertices contributes
with a real number between 0 and 4.
Subgraphs which are not matchable of course contribute with 0
and two disjoint edges as well as four points on a line contribute with 1.
The interesting part are the squares, which admit two matchings,
say $\mu_1$ and $\mu_2$:
then $\delta(\aa + \kk_{\cH}, \cH) = |1 - \eta(\aa(\mu_1 - \mu_2))|^2$;
in general, this may be any number between 0 and 2 but
if $\aa \in H^1(\cG,\Zz)$ then this is 0 or 4.

In order to prove this Theorem, we need an auxiliary result in linear algebra.
The proof of Lemma 2.3 (actually a rather straightforward computation)
may be found in [??].

{\nobf Lemma 2.3: }
{\sl Let $P(t) = t^n + a_1 t^{n-1} + \cdots + a_{n-1} t + a_n$
be the singular polynomial of $A$
(where $A$ is an arbitrary $n \times n'$ complex matrix).
Then
$$a_m = (-1)^m \sum_B |\det B|^2$$
where $B$ ranges over all $m \times m$ submatrices of $A$.}

{\nobf Proof of Theorem 2.2: }
Since balanced subgraphs of $\cG$ with $2m$ elements
correspond to $m \times m$ submatrices of $A$,
this is a consequence of Lemma 1.4 and Lemma 2.3.
\qed

We now describe another, more graphical, interpretation for Theorem 2.2.
We define a {\sl pipe system} of $\cG$ as an oriented pair
$\nu = (\mu_1,\mu_2)$ of matchings of a subgraph $\cH$ of $\cG$;
we call $\cH$ the {\sl support} of the pipe system.
Figure 2.1 shows an example of a pipe system:
we draw the edges of $\mu_2$ oriented from black to white
and the edges of $\mu_1$ from white to black
(unused edges are represented by dotted lines).
A pipe system is thus a collection of pipes (i.e., oriented edges of $\cG$)
such that, at each vertex,
there is either one pipe coming in and one pipe going out
or no pipe coming in and no pipe going out
(the water that comes in must go out
and you can not pipe too much water through a vertex).
We define the {\sl size} $|\nu|$ of the pipe system
as half the number of vertices in $\cH$
(in Figure 3.3, $|\nu| = 5$).

The Kasteleyn class $\kk_{\cH}$ shall be called $\kk_\nu$.
A pipe system obtains an element of $C_1(\cH,\ZZ)$
(and thus of $C_1(\cG,\ZZ)$)
but must not be confused with it:
if two pipes cancel each other homologically,
they still have to be taken into account for the pipe system.
If $\aa \in H^1(\cH,\SS^1)$,
$\eta(\aa(\nu))$ and $\eta((\aa + \kk_\nu)(\nu))$ are well defined
complex numbers.

\midinsert \smallskip

\centerline{\psfig{file=pipes.eps}}

\smallskip

\centerline{\eightrm Figure 2.2}

\smallskip \endinsert

{\nobf Corollary 2.4: }
{\sl Let $A$ be a generalized Kasteleyn matrix and let
$P_\aa(t) = t^n + a_1 t^{n-1} + \cdots + a_{n-1} t + a_n$
be the singular polynomial of $A$.
Then
$$P_\aa(t) = \sum_\nu (-1)^{|\nu|} t^{n-|\nu|} \eta((\aa + \kk_\nu)(\nu)),$$
where $\nu$ ranges over all pipe systems of $\cG$.}

{\nobf Proof: }
This follows directly from Theorem 2.2 and the definitions.
\qed

\bigbreak

{\nobf 3. Singular polynomials of planar graphs}

\nobreak\smallskip\nobreak

Theorem 2.2 provides an interpretation for the coefficients of
singular polynomials of arbitrary generalized Kasteleyn matrices.
In this section we take a closer look at the right hand side of equation (2)
when $A$ is a Kasteleyn matrix.

For planar balanced bipartite graphs $\cH \subseteq \cG$,
define $\pp_{\cG,\cH} \in H^1(\cH, \Zz)$
by $\pp_{\cG,\cH} = \pi_{\cG,\cH}(\kk_\cG) - \kk_\cH$.
In Figure 3.1 we illustrate the several objects involved in this definition:
(a), (b), (c) and (d) represent $\kk_\cG$,
$\pi_{\cG,\cH}(\kk_\cG)$, $\kk_\cH$ and $\pp_{\cG,\cH}$, respectively,
where again solid lines stand for a label 1 and
dashed lines stand for a label $-1$.
The following lemma provides an alternate definition for this class.

{\nobf Lemma 3.1: }
{\sl Let $\cH \subset \cG$ be balanced planar graphs
and let $C$ be a cycle in $\cH$.
Let $q$ be the number of vertices of $\cG$ not belonging to $\cH$
which are inside $C$. Then 
$$\pp_{\cG,\cH}(C) \equiv q \pmod 2.$$
}

{\nobf Proof: }
This follows directly from equation $(1)$ in Lemma 1.1.
\qed

\midinsert \smallskip

\centerline{%
\vbox{\hsize=25mm%
\centerline{\psfig{file=diag1a.eps}}
\bigskip \centerline{\eightrm (a)}}\qquad%
\vbox{\hsize=25mm%
\centerline{\psfig{file=diag1b.eps}}
\bigskip \centerline{\eightrm (b)}}\qquad%
\vbox{\hsize=25mm%
\centerline{\psfig{file=diag1c.eps}}
\bigskip \centerline{\eightrm (c)}}\qquad%
\vbox{\hsize=25mm%
\centerline{\psfig{file=diag1d.eps}}
\bigskip \centerline{\eightrm (d)}}}

\smallskip

\centerline{\eightrm Figure 3.1}

\smallskip \endinsert

In the hope of making the intuitive meaning of this definition clearer,
especially for adjacency graphs of quadriculated or triangulated disks,
we introduce some extra structure.

Let $\bcG$ be the CW-complex obtained from $\cG$
by closing each hole with a 2-cell;
$\bcG$ is thus always homeomorphic to a disk.
For $\cH \subseteq \cG$, let $\bcH \subseteq \bcG$
be the obtained from $\cH$ by adding the 2-cells of $\bcG$
whose boundaries are contained in $\cH$;
in other words, we close the holes of $\cH$ which contain no points of $\cG$.
The inclusion $\cH \subseteq \bcH$ induces
an injective map from $H^1(\bcH, \Zz)$ to $H^1(\cH, \Zz)$
which allows for a natural identification
of $H^1(\bcH, \Zz)$ with a subset of $H^1(\cH, \Zz)$.

{\nobf Lemma 3.2: }
{\sl $\pp_{\cG,\cH}$ belongs to $H^1(\bcH, \Zz)$.}

{\nobf Proof: }
This follows easily from Lemma 3.1.
\qed

Recall that two tilings by dominoes of a quadriculated region are said
to differ by a {\sl flip} if they coincide except for two dominoes;
in other words, their difference (in the homological sense) is a square.
If $\cG$ is the graph of a quadriculated planar region,
the difference between two tilings of $\cH$ differing by a flip
is 0 in $H_1(\bcH, \ZZ)$;
thus, tilings mutually accessible by flips always have the same 
$\pp_{\cG,\cH}$-parity.

We may now state the promised interpretation
for the coefficients of singular polynomial $P_\kk$ of $K$.
Since $P_\kk$ is well defined from $\cG$, we may adopt a lighter notation
and call it $P_\cG$, the {\sl singular polynomial} of $\cG$.

{\nobf Theorem 3.3: }
{\sl Let $\cG$ be a planar bipartite graph
and let $P_\cG = t^n + k_1 t^{n-1} + \cdots + k_{n-1} t + k_n$
be the singular polynomial of $\cG$. Then
$$k_m = (-1)^m \sum_{|\cH| = 2m}%
\delta(\pp_{\cG,\cH}, \cH) \eqno{(3)}$$
where $\cH$ ranges over all balanced subgraphs with $2m$ vertices.}

{\nobf Proof: }
This is a corollary of Theorem 2.2 and the definition of $\pp_{\cG,\cH}$.
\qed

Recall that if $\pp_{\cG,\cH} = 0$ then $\delta(\pp_{\cG,\cH}, \cH)$
is just the square of the number of matchings of $\cH$.
This always happens if $\bcH$ is simply connected.
As a corollary, if $\cG$ is the graph of a quadriculated planar region
and $m \le 3$,
or if $\cG$ is the graph of a triangulated planar region and $m \le 5$, then
$$k_m = (-1)^m \sum_{|\cH| = 2m} \delta(0, \cH)$$
where $\cH$ ranges over all balanced subgraphs (subregions)
with $2m$ vertices (squares, triangles).

Notice that $\pp_{\cG,\cH}$,
and thus the right hand side of equation (3),
depends on the way $\cG$ is drawn in the plane.
Examples (a) and (c) in Figure 1.1 show that $P_\cG$
indeed depends on the way $\cG$ is drawn:
for (a) we have $k_2 = 9$ but for (b) we have $k_2 = 5$.
This causes the singular values to change in a complicated way:
for (a) the singular values are approximately
$0.5549581321, 0.8019377358, 2.246979604$ while for (c) they are
$0.3472963553, 1.532088886, 1.879385242$.

It is natural to conjecture that the number of non-zero singular values
coincides with the size of a maximal partial matching of $\cG$.
In Figure 3.2(a) we present an example to show that this is not always true:
there are partial matchings of size 3 but since the singular
polynomial of $\cG$ is $t^3 - 7 t^2 + 10 t$ there are only two
non-zero singular values: $\sqrt 2$ and $\sqrt 5$.
In Figure 3.2(b) we draw the same graph in a different way
and we now have three non-zero singular values: 1, $\sqrt 2$ and 2.

\midinsert \smallskip

\centerline{%
\vbox{\hsize=45mm%
\centerline{\psfig{file=diag3a.eps}}
\bigskip \centerline{\eightrm (a)}}\qquad%
\vbox{\hsize=45mm%
\centerline{\psfig{file=diag3b.eps}}
\bigskip \centerline{\eightrm (b)}}}

\smallskip

\centerline{\eightrm Figure 3.2}

\smallskip \endinsert

We state Theorem 3.3 in the language of pipe systems.
Denote $\pp_{\cG,\cH}(\nu)$ (where $\cH$ is the support of $\nu$)
by $\pp(\nu) \in \Zz$.
We describe an elementary definition of $\pp(\nu)$.
Join pairs of vertices not in $\cH$,
matching black vertices with white vertices in an arbitrary way;
if all intersections are transversal,
$\pp(\nu)$ is the parity of the number of intersections
of such new lines with the pipes.

{\nobf Corollary 3.4: }
{\sl Let $\cG$ be a planar graph and $P_\cG(t)$ its singular polynomial.
Then
$$P_\cG(t) = \sum_\nu (-1)^{|\nu|} t^{n-|\nu|} \eta(\pp(\nu)),$$
where $\nu$ ranges over all pipe systems of $\cG$.}

{\nobf Proof: }
This follows directly from Theorem 3.3 and Corollary 2.4.
\qed

\bigbreak

{\nobf 4. $q$}

\nobreak

In several branches of combinatorics,
$q$-analogues or quantizations of classical problems
have been seen to be interesting and useful.
There are often several interpretations for the $q$-analogue
of a given concept, some sophisticated
(involving quantum groups and the like) and some elementary.
In this section we briefly consider a $q$-analogue of Kasteleyn matrices
in a very na{\"\i}ve way and extend the results of the previous sections
to this setting;
our interest in doing so is that the methods of the previous sections
extend very easily to this more general context and the coefficients
will actually have a natural interpretation.

Let $\Cq = \CC[q,q^{-1}]$.
We extend the usual complex conjugation to $\Cq$
by postulating $\bar q = q^{-1}$;
$q$ may be thought of as an unknown complex number of absolute value 1.
Let $\cG$ be a planar bipartite graph
with $n$ white vertices and $n'$ black vertices.
A {\sl generalized Kasteleyn $q$-matrix} for $\cG$ is an $n \times n'$
matrix $A$ with coefficients in $\Cq$ such that
$$a_{ij}\bar a_{ij} = \cases{1,&if the $i$-th white vertex
and the $j$-th black vertex are adjacent,\cr
0& otherwise;\cr}$$
thus, the entries are always monomials
and substituting $q$ by a complex number of absolute value 1
changes a generalized Kasteleyn $q$-matrix into an ordinary
generalized Kasteleyn matrix.

Consider the additive group $\Sq = \SS^1 \oplus \ZZ$ and let $\qq$
be the canonical generator of the $\ZZ$ component.
If we extend the classical $\eta$ to $\eta: \SS^1 \oplus \ZZ \to \Cq$
by postulating $\eta(\qq) = q$ we may identify
a generalized Kasteleyn $q$-matrix $A$ with a cocomplex $\AA \in C^1(\cG,\Sq)$.
Again, $C^2(\cG,\Sq) = 0$ and $A$ defines an element $\aa \in H^1(\cG,\Sq)$.

Let $\cG$ be a planar graph.
Bounded connected components of the complement of $\cG$
have well defined positively oriented boundaries $\beta$ in $H_1(\cG,\ZZ)$.
We define a Kasteleyn $q$-matrix to be
a generalized Kasteleyn $q$-matrix $A$ such that $\aa(\beta) = \qq$
for all such boundaries $\beta$.
We define the {\sl singular $q$-polynomial} of $\cG$
to be the singular polynomial of a Kasteleyn $q$-matrix of $\cG$.
As before, singular $q$-polynomials are easuly seen
to be well defined but now they are of course
polynomials in $\Cq[t]$, or, rather equivalently,
polynomials in two variables $q$ and $t$.
Finally, define the {\sl area} of a pipe system $\nu$, $A(\nu)$
to be the number of bounded connected components of the complement of $\cG$
positively surrounded by $\nu$,
counted with sign and multiplicity.

With these definitions we have the following theorem.

{\nobf Theorem 4.1: }
{\sl Let $\cG$ be a planar graph and $P_\cG(q,t)$ its singular $q$-polynomial.
Then
$$P_\cG(q,t) = \sum_\nu (-1)^{|\nu|} q^{A(\nu)} t^{n-|\nu|} \eta(\pp(\nu)),$$
where $\nu$ ranges over all pipe systems of $\cG$.}

Since the proof is entirely analogous to that of Corollary 3.4,
we leave the details to the reader.

\bigbreak

{\nobf 5. Rectangles}

\nobreak

Kasteleyn ([K]) computes the determinant of $K$
for rectangles essentially by computing its singular values.
For the reader's convenience,
we repeat that part of his work in our language.
In order to simplify notation in the statement and proof, let
$$\eqalign{
X^+_{M,N} &= \left\{ { (k,\ell) \in \ZZ^2 \;\big|\;
\hbox{$\displaystyle 1 \le k \le {M+1\over 2}$ and if
$\displaystyle k = {M+1\over 2}$ then
$\displaystyle 1 \le \ell \le {N+1\over 2}$}
} \right\},\cr
X^-_{M,N} &= \left\{ { (k,\ell) \in \ZZ^2 \;\big|\;
\hbox{$\displaystyle 1 \le k \le {M+1\over 2}$ and if
$\displaystyle k = {M+1\over 2}$ then
$\displaystyle 1 \le \ell < {N+1\over 2}$}
} \right\}.\cr}$$

{\nobf Theorem 5.1: }
{\sl Let $\cG$ be a $M \times N$ rectangular grid
and let $K$ be its Kasteleyn matrix.
Then the non-zero singular values of $K$ are $\sigma_{k,\ell}$,
$(k,\ell) \in X^-_{M,N}$, where
$$\sigma_{k,\ell}^2 =
(\alpha^k + \alpha^{-k})^2 + (\beta^\ell + \beta^{-\ell})^2,\quad
\alpha = \exp\left({\pi i\over M+1}\right),\quad
\beta = \exp\left({\pi i\over N+1}\right).$$
}

The complicated description of the allowed values of the indices $k$ and $\ell$
is necessary in order to avoid zeroes and duplications
in a way which is correct for all possible parities of $M$ and $N$
(Kasteleyn has a simpler formula since he assumes $N$ to be even).
Notice that
$$\sigma_{k,\ell} =
2\left({\cos^2{k\pi \over M+1} + \cos^2{\ell\pi \over N+1}}\right)^{1/2},$$
(an expression closer to Kasteleyn's),
$\sigma_{k,\ell} = \sigma_{M+1-k,\ell} =
\sigma_{k,N+1-\ell} = \sigma_{M+1-k,N+1-\ell}$
and that $\sigma_{k,\ell} = 0$ if and only if $M$ and $N$ are both odd,
$k = (M+1)/2$ and $\ell = (N+1)/2$.
Thus, if we just demand $1 \le k \le M$ and $1 \le \ell \le N$
then all non-zero singular values are counted twice
and we occasionally introduce a 0.

{\nobf Proof: }
We index vertices by pairs $(k',\ell')$, $1 \le k' \le m$, $1 \le \ell' \le n$.
The vertex $(k',\ell')$ is called white when $k' + \ell'$ is even.
Define $K$ as the Kasteleyn matrix with entries 1 for horizontal edges
and $i$ for vertical edges: $K$ defines a linear transformation
from the ``black space'' to the ``white space''.
Consider the white vectors
$$w_{k,l} = (\alpha^{kk'} - \alpha^{-kk'})%
(\beta^{\ell\ell'} - \beta^{-\ell\ell'}),$$
$(k,\ell) \in X^+_{m,n}$:
they clearly form an orthogonal basis for the white space
(this is where a careful choice of $X^+_{m,n}$ becomes necessary).
Similarly, the black vectors $b_{k,l}$ defined by the same formula with
$(k,\ell) \in X^-_{m,n}$ form an orthogonal basis for the black space.
A simple computation yields
$|w_{k,l}| = |b_{k,l}|$ for $(k,l) \in X^-_{m,n}$ and
$$\eqalign{
K b_{k,\ell} &=
\left({(\alpha^k + \alpha^{-k}) + i (\beta^\ell + \beta^{-\ell})}\right)
w_{k,\ell},\cr
K^\ast w_{k,\ell} &=
\left({(\alpha^k + \alpha^{-k}) - i (\beta^\ell + \beta^{-\ell})}\right)
b_{k,\ell}.\cr
}$$
Thus, $w_{k,\ell}$ and $b_{k,\ell}$ are singular vectors
and $\sigma_{k,\ell}$ are singular values.
\qed

{\nobf Corollary 5.2: }
{\sl Let $\cG$ be a $m \times n$ rectangular grid. Let
$$\alpha = \exp\left({\pi i\over m+1}\right),\quad
\beta = \exp\left({\pi i\over n+1}\right),\quad
N = \left\lfloor{mn \over 2}\right\rfloor.$$
Then
$$\prod_{(k,\ell) \in X^-_{m,n}}%
(t - (\alpha^k + \alpha^{-k})^2 - (\beta^\ell + \beta^{-\ell})^2) =
\sum_{j = 0 \ldots N}t^{N-j} (-1)^j \sum_{|\cH| = 2j}%
\delta(\pp_{\cG,\cH}, \cH),$$
where $\cH$ ranges over all balanced subgraphs with $2j$ vertices.}

{\nobf Proof: }
This follows directly from Theorem 3.3 and Theorem 4.1.
\qed

These results show that the characteristic polynomial of $KK^\ast$
usually factors a lot if $\cG$ is a rectangle.
If $\zeta$ is a root of unity whose order $M$
is the least common multiple of $2(m+1)$ and $2(n+1)$ then
all the roots $\sigma_{k,\ell}^2$ of this polynomial
are in $\RR \cap \ZZ[\zeta]$, a ring of degree $\phi(M)/2$ over $\ZZ$.
In particular, for square grids of order $n$,
irreducible factors of the characteristic polynomial of $KK^\ast$
have degree at most $n$.
Here are a few sample examples; we give the polynomial
$\det(t I - KK^\ast) = t^{n} + k_1 t^{n - 1} + \cdots + k_{n-1} t + k_n$
(whose roots are the squares of singular values)
factored in $\ZZ$.

$$1 \leqno{ [1,1]}$$

$$t-1 \leqno{ [2,1]}$$

$$\left (t-2\right )^{2} \leqno{ [2,2]}$$

$$t-2 \leqno{ [3,1]}$$

$$\left (t-1\right )\left (t-3\right )^{2} \leqno{ [3,2]}$$

$$\left (t-2\right )^{2}\left (t-4\right )^{2} \leqno{ [3,3]}$$

$${t}^{2}-3\,t+1 \leqno{ [4,1]}$$

$$\left ({t}^{2}-5\,t+5\right )^{2} \leqno{ [4,2]}$$

$$\left ({t}^{2}-3\,t+1\right )\left ({t}^{2}-7\,t+11\right )^{2}
\leqno{ [4,3]}$$

$$\left ({t}^{2}-6\,t+4\right )^{2}\left (t-3\right )^{4} \leqno{ [4,4]}$$

$$\left (-1+t\right )\left (t-3\right ) \leqno{ [5,1]}$$

$$\left (t-1\right )\left (t-2\right )^{2}\left (t-4\right )^{2}
\leqno{ [5,2]}$$

$$\left (t-1\right )\left (t-2\right )\left (t-5\right )^{2}\left (t-3
\right )^{3}
\leqno{ [5,3]}$$

$$\left ({t}^{2}-3\,t+1\right )\left ({t}^{2}-5\,t+5\right )^{2}
\left ({t}^{2}-9\,t+19\right )^{2}
\leqno{ [5,4]}$$

$$\left (t-1\right )^{2}\left (t-2\right )^{2}\left (t-3\right )^{2}
\left (t-6\right )^{2}\left (t-4\right )^{4}
\leqno{ [5,5]}$$

$${t}^{3}-5\,{t}^{2}+6\,t-1 \leqno{ [6,1]}$$

$$\left ({t}^{3}-8\,{t}^{2}+19\,t-13\right )^{2} \leqno{ [6,2]}$$

$$\left ({t}^{3}-5\,{t}^{2}+6\,t-1\right )\left ({t}^{3}-11\,{t}^{2}+38\,
t-41\right )^{2}
\leqno{ [6,3]}$$

$$\left ({t}^{6}-19\,{t}^{5}+142\,{t}^{4}-529\,{t}^{3}+1017\,{t}^{2}-922
\,t+281\right )^{2}
\leqno{ [6,4]}$$

$$\left ({t}^{3}-5\,{t}^{2}+6\,t-1\right )\left ({t}^{3}-8\,{t}^{2}+19\,t
-13\right )^{2}\left ({t}^{3}-14\,{t}^{2}+63\,t-91\right )^{2}
\leqno{ [6,5]}$$

$$\left ({t}^{3}-10\,{t}^{2}+24\,t-8\right )^{2}\left ({t}^{3}-10\,{t}^{2
}+31\,t-29\right )^{4}
\leqno{ [6,6]}$$

$$\left (t-2\right )^{2}\left ({t}^{2}-4\,t+2\right )^{2}\left ({t}^{2}-8
\,t+8\right )^{2}\left ({t}^{2}-8\,t+14\right )^{4}\left (t-4\right )^{
6}
\leqno{ [7,7]}$$

$$(t - 2)^2 (t^3 - 12 t^2 + 36 t - 8)^2
(t^3 - 9 t^2 + 24 t - 17)^4 (t^3 - 12 t^2 + 45 t - 53)^4
\leqno{ [8,8]}$$

\bigbreak

{\nobf 6. Examples}

\nobreak

Although Aztec diamonds have so many interesting properties
(see [EKLP] and [P]),
the characteristic polynomial of $KK^\ast$ does not factor very much:

$$
(t^3 - 8 t^2 + 17 t - 8)^2
\leqno\hbox{2-Aztec diamond}
$$

$$
(t^4 - 10 t^3 + 28 t^2 - 24 t + 4)^2 (t - 4)^4
\leqno\hbox{3-Aztec diamond}
$$

$$
(t^{10} - 32 t^{9} + 441 t^{8} - 3424 t^{7} + 16432 t^{6} - 50240 t^{5}
+ 97041 t^{4} - 112896 t^{3} + 70921 t^{2} - 18784 t + 1024)^{2}
\leqno\hbox{4-Aztec diamond}
$$

$$
(t^{11} - 34 t^{10} + 496 t^{9} -4064 t^{8} +20562 t^{7}
- 66524 t^{6} + 137728 t^{5} 
\atop
- 177120 t^{4} + 131825 t^{3} - 49066 t^{2} + 6576 t - 128)^2
(t - 4)^8
\leqno\hbox{5-Aztec diamond}
$$

The fact that these polynomials are always squares follows from symmetry.
The factor $(t - 4)^{4k}$ seems to appear in the $2k+1$-Aztec diamond,
a fact for which we have no explanation.

Finally, here is a small list of ``irregular'' examples.
A possible real Kasteleyn matrix is indicated by the dashed lines (the $-1$'s).

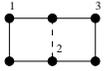
$$
t^3 - 7 t^2 + 15 t - 9 = (t - 1)(t - 3)^2
\atop 1.732051, 1.732051, 1 
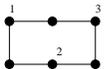
\leqno{\psfig{file=diag2e.eps}}
$$

$$
t^3 - 6 t^2 + 9 t - 4 = (t-4)(t-1)^2 \atop 2, 1, 1
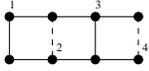
\leqno{\psfig{file=diag2f.eps}}
$$

$$
t^4 - 10 t^3 + 35 t^2 - 50 t + 25 = (t^2 - 5 t + 5)^2
\atop 
1.902113, 1.902113, 1.175571, 1.175571
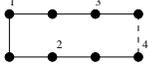
\leqno{\psfig{file=diag2g.eps}}
$$

$$
t^4 - 8 t^3 + 20 t^2 - 16 t + 4 = (t^2 - 4 t + 2)^2\atop
1.847759, 1.847759, 0.765367, 0.765367
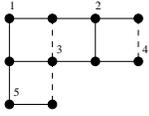
\leqno{\psfig{file=diag2h.eps}}
$$

$$
\displaystyle {t^{5} - 13 t^4 + 63 t^3 - 140 t^2 + 140 t - 49 = \atop
(t^2 - 6 t  + 7) (t^3 - 7 t^2  + 14 t  - 7)}
\atop
2.101003, 1.949856, 1.563663, 1.259280, 0.867768
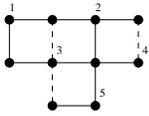
\leqno{\psfig{file=diag2i.eps}}
$$

$$
\displaystyle {t^{5} - 13 t^4 + 62 t^3 - 132 t^2 + 121 t - 36 = \atop
(t - 4)(t^4 - 9 t^3 + 26 t^2 - 28 t + 9)}
\atop
2.126757, 2, 1.576415, 1.197126, 0.747468
\leqno{\psfig{file=diag2j.eps}}
$$
\bigbreak\vfil

{\noindent \bigbf References:}

\nobreak
\smallskip
\nobreak

\parindent=40pt

\item{[EKLP]}{N. Elkies, G. Kuperberg, M. Larsen and J. Propp,
{Alternating-sign matrices and domino tilings},
Journal of Algebraic Combinatorics, {\bf 1}, 111-132 and 219--234 (1992).}

\item{[K]}{P. W. Kasteleyn,
{The statistics of dimers on a lattice I. The number of dimer
arrangements on a quadratic lattice,}
Phisica {\bf 27}, 1209-1225 (1961).}
 
\item{[LL]}{E. H. Lieb and M. Loss,
{Fluxes, Laplacians and Kasteleyn's theorem,}
Duke Math. Jour., {\bf 71}, 337-363 (1993).}

\item{[P]}{J. Propp,
{Twenty open problems in enumeration of matchings,}
preprint (1996).}

\item{[ST]}{N. C. Saldanha and C. Tomei,
{An overview of domino and lozenge tilings,}
Resenhas IME-USP, {\bf2}(2), 239--252, (1995).}
	
\bigskip

\vbox{
\obeylines
\parskip=0pt
\parindent=0pt
Nicolau C. Saldanha
\smallskip
Departamento de Matem\'atica, PUC-Rio
Rua Marqu\^es de S\~ao Vicente 225
G\'avea, Rio de Janeiro, RJ 22453-900, BRAZIL
\smallskip
UMPA, ENS-Lyon, 46 Allee d'Italie
69364 Lyon cedex 07, FRANCE
\smallskip
nicolau@mat.puc-rio.br
http://www.mat.puc-rio.br/$\sim$nicolau/
}

\bye